\documentclass[12pt,leqno]{amsart}
\usepackage{amscd}
\usepackage{amssymb}
\usepackage[matrix,arrow]{xy}
\usepackage{geometry}
\usepackage{color}

\newtheorem{theorem}{Theorem}[section]
\newtheorem{proposition}[theorem]{Proposition}
\newtheorem{lemma}[theorem]{Lemma}

\newtheorem{corollary}[theorem]{Corollary}
\theoremstyle{remark}
\newtheorem{remark}[theorem]{Remark}

\newtheorem{definition}[theorem]{Definition}

\newcommand{\uk}{\mathrm{\underline{K}}}

\newcommand{\ot}{\otimes}

\newcommand{\CCC}{\mathbb{C}}
\newcommand{\AAA}{\mathcal{E}}

\newcommand{\LLL}{\mathcal{L}}
\newcommand{\ep}{\varepsilon}
\newcommand{\st}{such that}

\newcommand{\BBB}{\mathcal{B}}

\newcommand{\OOO}{\mathcal{O}}

\newcommand{\Ad}[1]{\operatorname{Ad}({#1})}

\newcommand{\fg}{finitely generated}

\newcommand{\uka}{\underline{K}(A)}
\newcommand{\ukb}{\underline{K}(B)}
\newcommand{\fset}{\mathcal{F}}

\newcommand{\pset}{\mathcal{P}}

\newcommand{\gset}{\mathcal{G}}

\newcommand{\kalg}{Kirchberg algebra}

\newcommand{\kalgs}{Kirchberg algebras}

\newcommand{\inv}{\mathrm{Inv}}
\newcounter{rocount}

\begin{document}
\title{Continuous Fields of Kirchberg C*-algebras}
\dedicatory{Dedicated to George Elliott on his 60th birthday}
\author{Marius Dadarlat and Cornel Pasnicu}
\address{Department of Mathematics, Purdue University, West
Lafayette IN 47907, U.S.A.} \email{mdd@math.purdue.edu}
\address{University of Puerto Rico, Department of Mathematics, Box 23355, San
Juan, PR 00931, U.S.A.} \email{cpasnic@upracd.upr.clu.edu}
\date{\today}

\begin{abstract} In this paper we study the C*-algebras associated to
continuous fields over locally compact metrisable zero dimensional
spaces whose fibers are Kirchberg C*-algebras satisfying the UCT.
We show that these algebras are inductive limits of finite direct
sums of Kirchberg algebras and they are classified up to
isomorphism by topological invariants.
\end{abstract}

\maketitle

\section{Introduction}
A  purely infinite separable simple  nuclear C*-algebra is called
a Kirchberg algebra. Kirchberg \cite{Kir:class} and Phillips
\cite{Phi:class} proved that two Kirchberg algebras $A$ and $B$
are stably isomorphic if and only if they are KK-equivalent.
Consequently, if in addition $A$ and $B$ satisfy the universal
coefficient theorem (UCT) of \cite{RosSho:UCT}, then $A \cong B$
if and only if $K_*(A)\cong K_*(B)$.

Kirchberg \cite{Kir:Michael} generalized the KK-theory isomorphism
result to  nonsimple C*-algebras. He showed that if  $A$ and $B$
are nuclear separable C*-algebras with primitive ideal spectrum
homeomorphic to some $\mathrm{T}_0$-space $X$, then
\mbox{$A\otimes\mathcal{O}_\infty \otimes\mathcal{K}\cong
B\otimes\mathcal{O}_\infty \otimes\mathcal{K}$} if and only if $A$
is $KK_X$-equivalent to $B$, where $KK_X$ is a suitable
generalization of the Kasparov theory which preserves
 the primitive ideal spectra of $A$ and
$B$ (or rather the lattice homomorphisms from the open subsets of
$X$ to the lattice of closed ideals of $A$ and $B$). Unlike the
case of simple C*-algebras, as observed in \cite{KR2}, one does
not have a general algebraic criterion for recognizing when two
C*-algebras are $KK_X$-equivalent, as we lack a generalization of
the UCT for $KK_X$. Finding such a criterion seems to be a
difficult problem even when the space $X$ consists of finitely
many points and is non Hausdorff. The case when $X$ consists of
two points  was solved by R{\o}rdam in \cite{Ror:6-term}.

One of the goals of  the present paper is to propose an answer to
the above question for the separable nuclear  C*-algebras $A$
whose primitive ideal spectrum $\mathrm{Prim}(A)$ is  zero
dimensional and Hausdorff, under the assumption that all simple
quotients of $A$ satisfy the UCT. We introduce a homotopy
invariant $\inv(A)$ consisting of a preordered semigroup $$P(A
\otimes \mathcal{O}_2) \oplus \underline{K} (A),$$ together with
the action of the Bockstein operations on $\uka$.
 Here $P(A \otimes \mathcal{O}_2)$
denotes the Murray-von Neumann semigroup of equivalence classes of
projections in $A \otimes \mathcal{O}_2\otimes \mathcal{K}$ and
$\uka$ is the total K-theory group of $A$, see Section~\ref{invar}.
It turns out that if $A$ and $B$ are separable nuclear C*-algebras
with zero dimensional Hausdorff primitive spectra and with all
simple quotients satisfying the UCT, then
\mbox{$A\otimes\mathcal{O}_\infty \otimes\mathcal{K}\cong
B\otimes\mathcal{O}_\infty \otimes\mathcal{K}$} if and only
\mbox{$\inv(A)\cong \inv(B)$}. Since we do not use Kirchberg's
$KK_X$-theory, we deduce that for such algebras, $A$ is
$KK_X$-equivalent to $B$ if and only if $\inv(A)\cong \inv(B)$. The
examples  constructed in \cite{DadEil:bock} show that the action of
the Bockstein operations is an essential part of the invariant.
$\inv(A)$ is an adaptation to purely infinite C*-algebras of an
invariant introduced in \cite{DadLor:duke} and
\cite{DadGong:class-rr0}. A positive morphism $\inv(A)\to \inv(B)$
must preserve the filtration of the total K-theory group $\uk(-)$
induced by ideals. Since all the ideals of the C*-algebras
classified in this paper give rise to quasidiagonal extensions, the
corresponding boundary maps vanish, and hence there is no need to
include them in the classifying invariant as it was necessary to do
in \cite{Ror:6-term}.

 A related goal of this paper is to describe the structure and the
classification of the C*-algebras associated to continuous fields
over locally compact metrisable zero dimensional spaces whose
fibers are Kirchberg algebras satisfying the UCT. We show that
these algebras are inductive limits of finite direct sums of
Kirchberg algebras satisfying the UCT and that they are classified
up to isomorphism by the topological invariant $\inv_u(-): = (\inv
(-), P_u(-), \tau)$, see Sections ~\ref{invar} and
~\ref{section-class}. It is worth to note that their structure is
obtained proving first the classification result for a larger
class of C*-algebras, by the same invariant. Since $\inv_u$ is
continuous and homotopy invariant, we deduce immediately that in
this (larger) class of C*-algebras, the isomorphism and the shape
equivalence are equivalent properties, see Section
~\ref{section-class}.

 Let us describe how the paper is organized.
Section~\ref{section-prelim} is devoted to preliminaries. Using
the results of \cite{Kir:class} and \cite{Phi:class} we show in
Section \ref{section-fields} that the C*-algebra associated to a
continuous field of Kirchberg algebras satisfying the UCT over a
metrisable zero dimensional locally compact space admits local
approximations by finite direct sums of Kirchberg algebras
satisfying the UCT and having finitely generated K-theory.
 In
Section~\ref{invar} we introduce the  invariants $\inv$ and
$\inv_u$ and describe their basic properties.
 In
Section~\ref{section-class} we prove that the C*-algebras $A$
which admit local approximations as described above are classified
up to stable isomorphism by $\inv(A)$ and in fact they can be
written as inductive limits of finite direct sums of Kirchberg
algebras satisfying the UCT. In particular, these structure and
classification results apply to the C*-algebras associated to
continuous fields of Kirchberg algebras satisfying the UCT over a
metrisable zero dimensional locally compact space.

\section{Preliminaries}\label{section-prelim}

\begin{definition}\label{2.3}
  A sequence $(A_n)$ of C*-subalgebras of a C*-algebra $A$ is
  called exhaustive if for any finite subset $\fset \subset A$, any $m>0$, and
  any $\ep>0$ one has $\fset\subset_{\ep}A_n$ for some $n>m$.
  The inclusion maps $A_n\hookrightarrow A$ are denoted by
  $\imath_n$ ( respectively $\jmath_n:B_n \hookrightarrow B$ if $(B_n)$
   is an exhausting sequence
  for $B$).
\end{definition}
Let $A$ be a separable C*-algebra and let $(A_n)$ be an exhaustive
sequence for $A$. Let $\{x_1,x_2,\dots,x_n,\dots\}\subset A$ be a
dense subset of $A$. After passing to a subsequence of $(A_n)$ we
may arrange that
\begin{equation}\label{exhaustive}
    \{x_1,x_2,\dots,x_k\}\subset_{1/k}A_k,\quad \forall\,k\in
\mathbb{N}.
\end{equation}
It the sequel we will always work with exhaustive sequences
satisfying \eqref{exhaustive}.

\begin{lemma}\label{5.5} Let $A$ be a separable $C^{*}$-algebra with an
exhaustive sequence $(A_{n})$ consisting of unital, nuclear
$C^{*}$-algebras.  Then after passing to a subsequence of
$(A_{n})$  satisfying \eqref{exhaustive} there is a sequence of
completely positive contractions $\mu_n : A \rightarrow A_n$ which
is asymptotically multiplicative and such that $lim_{n\rightarrow
\infty}\|\imath_{n}\mu_{n}(a)-a\| =0$ for all $a\in A$, where
$\imath_{n} : A_{n} \rightarrow A$ are the inclusion maps.
\end{lemma}
\begin{proof}Let $\mathcal{F}:=\{x_{1},x_{2},...,x_{n},...\}$ with
$\mathcal{\overline{F}}=A$, and let $\mathcal{F}_{k}:=\{ x_{1},
x_{2},...,x_{k} \} (k\in \mathbb{N})$. After passing to a
subsequence we may assume that $(A_{n})$ satisfies
\eqref{exhaustive}. Let $\mathcal{G}_{k}$ be a finite subset of
$A_{k}\, (k\in \mathbb{N})$ such that:
\begin{equation}\label{13}
\mathcal{F}_{k}\subseteq_{1/k} \mathcal{G}_{k}, \quad k\in
\mathbb{N}.
\end{equation}

For each $n$, let $e_n$ be the unit of $A_n$.  Using   that
$\mathcal{\overline{F}}=A$ and \eqref{13} it is easy to check that
$(e_{n})$ is a (\emph{not necessarily increasing}) approximate
unit of $A$. Since each $A_{n}$ is nuclear and unital, it follows
that we can find an approximate factorization of $id_{A_{n}}$ by
unital completely positive maps, on the finite set
$\mathcal{G}_{n}$ within ${1/n}$.
\[
\xymatrix{
{A_n}\ar@{=}[r]^{id}\ar[dr]_{f_n}&{A_n}\\
&{M_{k(n)}(\CCC)}\ar[u]_{g_n}}
\]

\begin{equation}\label{14}
    \|g_nf_n(y)-y\|<1/n,\quad y \in \gset_n,\, n\in \mathbb{N}.
\end{equation}
Now, using Arveson's extension theorem for unital completely
positive maps (see e.g. \cite[Thm. 6.1.5]{Ror:encyclopedia}), we
extend $f_n$ to a unital completely positive map (denoted in the
same way) $f_{n}:e_{n}Ae_{n}\rightarrow M_{k(n)}(\CCC)$ for every
$n\in \mathbb{N}$. Define  $h_n:A \to e_nAe_n$ by $h_n(a)=e_n a
e_n$ and $\mu_{n}:A\rightarrow A_{n}$ by $\mu_{n}=g_nf_nh_n$.
 Hence we have a diagram:
\[
\xymatrix{A\ar[rrd]^{\mu_n}\ar[ddr]_{h_n}\\
&{A_n}\ar@{^{(}->}[d]\ar@{=}[r]^{id}\ar[dr]_{f_n}&{A_n}\\
&{e_nAe_n}\ar[r]_{f_n}&{M_{k(n)}(\CCC)}\ar[u]_{g_n}}
\]
    Fix now an arbitrary $x\in A$ and
    $\varepsilon >0$.  Since $\overline{\mathcal{F}}=A$, it follows
    that there is a positive integer $k$ such that
    $\|x-x_{k}\| < \varepsilon$.  Since
    $x_{k} \in \mathcal{F}_{k} \subseteq \mathcal{F}_{n}$ for
    every $n\geq k$, it follows by \eqref{13} that there is
    an element $y_{n} \in \mathcal{G}_{n} \subset
    A_{n}$ such that $\| x_{k}-y_{n}
    \| < {1/n}$.  Hence:
    \begin{equation}\label{15}\| x-y_{n}\| < \varepsilon +
    {1/n}, \forall\, n\geq k.\end{equation}
    Since as noticed earlier $lim_{n\rightarrow \infty}\| e_{n}xe_{n}-x\|
    = 0$, there is $k_1 \geq k$
    such that:
    \begin{equation}\label{16}\| e_{n}xe_{n}-x\| < \varepsilon,
    \forall\,
    n\geq k_{1}.\end{equation}
    Now, using \eqref{14}, \eqref{15}, and \eqref{16} we can write for every $n\geq
    k_{1}:$
    $$\| \mu_{n}(x)-x\| < \|
    \mu_{n}(x)-y_{n}\| + \varepsilon + {1/n} =
    \| g_{n}f_{n}h_{n}(x)-y_{n}\| +
    \varepsilon + {1/n}$$
    $$=\| g_{n}f_{n}(e_{n}xe_{n})-y_{n}\| + \varepsilon +{1/n}$$
    $$<\| g_{n}f_{n}(x)-y_{n}\| +
    2\varepsilon + {1/n}$$
    $$<\|g_{n}f_{n}(y_{n})-y_{n}\| +
    3\varepsilon + 2/n<
    3\varepsilon + 3/n.$$
Hence $lim_{n\rightarrow \infty} \| \imath_{n}\mu_n(x)-x\| = 0$.
Note also   that each $\mu_{n}:A\rightarrow A_{n}$ is a completely
positive contraction  as a composition of completely positive
contractions. Finally, the fact that $(\mu_{n})$ is asymptotically
multiplicative follows easily since $\imath_n$ are
$*$-monomorphisms and $lim_{n\rightarrow
\infty}\|\imath_{n}\mu_{n}(a)-a\| =0$ for all $a \in A$.
\end{proof}

\begin{definition}\label{2.5}
 Let $\mathcal{A}=((A(x))_{x\in
X},\Gamma)$ be a continuous field of C*-algebras over a metrisable
locally compact space $X$. Here $\Gamma$ consists of vector fields
$a$ on $X$, (i.e. $a \in \prod_{x\in X}A(x)$, $a(x)\in A(x)$)
satisfying a number of natural axioms, including the continuity of
the map $x \mapsto \|a(x)\|$, and the condition that $\Gamma$ is
closed under local uniform approximation \cite[10.1.2]{Dix:C*}. If
$X$ is compact, then $\Gamma$ is a C*-algebra, called the
C*-algebra associated to $\mathcal{A}.$ If $X$ is just locally
compact, the C*-algebra associated to $\mathcal{A}$ consists of
those vector fields $a \in \Gamma$ with the property that the map
$x \mapsto \|a(x)\|$ is vanishing at infinity. Given a family of
C*-algebras $(A(x))_{x \in X}$ there are in general many choices
for  $\Gamma$ which makes $((A(x))_{x\in X},\Gamma)$ a continuous
field of C*-algebras with non-isomorphic associated C*-algebras.

Let $Y \subseteq X$ be a closed subspace and let $\Gamma_Y$ be the
restriction of $\Gamma$ to $Y$. One verifies that
$\mathcal{A}|_Y:=((A(x))_{x\in Y},\Gamma_Y)$ is a continuous field
of C*-algebras on $Y$ \cite[10.1.12]{Dix:C*}.  We shall denote by
$A(Y)$ the C*-algebra associated to $\mathcal{A}|_Y.$ Observe that
$A(X)$ is the C*-algebra associated to $\mathcal{A}$.
\end{definition}
\begin{remark}\label{2.6}
 Let $\mathcal{A}=((A(x))_{x\in
X},\Gamma)$ be a continuous field of C*-algebras over a metrisable
zero dimensional locally compact space $X$. If $Y\subseteq X$ is
closed, then $$I(Y)= \{f \in A(X): f|_{Y} = 0 \}$$ is a closed,
two-sided ideal of $A(X)$. Using \cite[Prop. 10.1.12]{Dix:C*} one
shows that $A(X)/I(Y)\cong A(Y)$.
 Let
$(F_n)_{n=1}^\infty$ be a decreasing sequence of compact subsets
of $X$ forming a basis of neighborhoods of a point $x_0\in X$.
Then $A(x_0)\cong \varinjlim (A(F_n), {\Phi_n})$, where each
$\Phi_n : A(F_n) \rightarrow A(F_{n+1})$ is the restriction map
$(\Phi_n (f)= f|_{F_{n+1}}, \, f \in A(F_n))$. Indeed, the
$I(F_n)$'s form an increasing sequence of closed, two-sided ideals
of $A(X)$ and $\overline{\bigcup_{n \in \mathbb{N}} I(F_n)} =
I(x_0),$ since $\bigcap_{n \in \mathbb{N}}F_n=\{x_0\}$. Therefore:
$$ A(x_0) \cong A(X)/{I(x_0)} \cong \varinjlim
(A(X)/{I(F_n)}, \varphi_n) \cong \varinjlim (A(F_n), \Phi_n) $$
where $\varphi_n : A(X)/{I(F_n)} \rightarrow A(X)/{I(F_{n+1})}$ is
induced by the inclusion $I(F_n)\subseteq I(F_{n+1})$.
\end{remark}
\begin{remark}\label{2.7}
Let $\mathcal{A} = ((A(x))_{x \in X}, \Gamma)$ be a continuous
field of separable C*-algebras with a countable approximate unit
of projections over a metrisable zero dimensional locally compact
space X.  Let $U$ be a closed subset of $X$. Then $A(U)$ has a
countable approximate unit of projections. This follows from
\cite{P2}.
If all the $A(x)$'s have real
rank zero, then $A(U)$ has real rank zero by \cite[Thm. 2.1]{P3}.
\end{remark}

\section{Continuous fields of Kirchberg
algebras}\label{section-fields} In this section we establish an
approximation property for C*-algebras associated to continuous
fields of Kirchberg C*-algebras over zero dimensional spaces, see
Thm. \ref{3.3}.
\begin{definition}\label{2.1}
A separable nuclear simple purely infinite C*-algebra  is called a
\kalg.
\end{definition}
We refer the reader to \cite{Ror:encyclopedia} for a background
discussion of Kirchberg algebras.
\begin{remark}\label{2.2}
Any \kalg\ is either unital or of the form $A_0\otimes
\mathcal{K}$ where $A_0$ is a unital \kalg\ and $[1_{A_0}]=0$ in
$K_0(A_0)$. In particular all Kirchberg algebras admit an
approximate unit consisting of projections.
\end{remark}
We introduce here notation for certain classes of C*-algebras.
This notation is used to shorten the statements of certain
intermediate results.
\begin{definition}\label{2.4}
$\AAA$ consists of unital \kalgs.

 $\AAA_{uct}$ consists of unital \kalgs\ satisfying the UCT.

 $\AAA_{fg}$ consists of  unital \kalgs\  $A$ with $K_*(A)$ \fg.

 $\AAA_{fg-uct}=\AAA_{fg}\cap\AAA_{uct}$

$\BBB$ consists of finite direct sums of unital \kalgs. One
defines similarly $\BBB_{uct}$, $\BBB_{fg}$ and  $\BBB_{fg-uct}$.
A (nuclear) separable C*-algebra is in the class $\LLL$
    if it admits an exhaustive sequence  $(A_n)$ with each $A_n$
    in $\BBB$. One defines
similarly $\LLL_{uct}$, $\LLL_{fg}$ and  $\LLL_{fg-uct}$.
\end{definition}

\begin{lemma}\label{3.1}
Let $\mathcal{A} = ((A(x))_{x \in X} , \Gamma)$ be a continuous
field of Kirchberg algebras over a metrisable zero dimensional
locally compact space $X$.  Let $A$ be the $C^*$-algebra
associated to $\mathcal{A}$.  Then $A \cong A\otimes
\mathcal{O}_\infty$.
\end{lemma}
\begin{proof}
If we set $ B:= A \otimes \mathcal{K},$ then $B$ is separable,
stable, nuclear, and  $\mathrm{Prim}(B) \cong\mathrm{Prim}(A)
\cong X$ is Hausdorff and zero dimensional. Every nonzero simple
quotient of $B$ is of the form $A(x_0) \otimes \mathcal{K}$, for
some $x_0 \in X$. But $A(x_0)$ being a Kirchberg algebra is purely
infinite and since the property of being purely infinite is
invariant under stable isomorphism \mbox{\cite[Thm. 4.23]{KR1}} it
follows that $ A(x_0) \otimes \mathcal{K}$ is also purely
infinite. Now, by \cite[Thm. 1.5]{BK} we deduce that $B$ is
strongly purely infinite. Using again the fact that $A$ and $B$
are stably isomorphic, \cite[Prop. 5.11 (iii)]{KR2} implies that
$A$ is strongly purely infinite. Hence the $C^*$-algebra $A$ is
separable nuclear, has an approximate unit of projections (see
Remark~\ref{2.7}) and is strongly purely infinite. Then, by
\cite[Thm. 8.6]{KR2} we have that $A \cong A \otimes
\mathcal{O}_\infty .$
\end{proof}
\begin{lemma}\label{3.2}
Let $\mathcal{A}=((A(x))_{x \in X} , \Gamma)$ be a continuous
field of unital Kirchberg algebras satisfying the $UCT$ over a
metrisable zero dimensional locally compact space.  Let $A$ be the
$C^*$-algebra associated to $\mathcal{A}$.  Let $\mathcal{F}
\subset A$ be a finite set and let $\varepsilon > 0$.  Then, for
every $x \in X$, there is a clopen neighborhood $U$ of $x$, and
there is a unital Kirchberg algebra $B$ satisfying the $UCT$, with
$K_* (B)$ finitely generated and a $*$-homomorphism $\gamma : B
\rightarrow A(U)$ such that \mbox{dist $(a|_U, \gamma (B)) <
\varepsilon$} for all $a \in \mathcal{F}$.  The distance is
calculated in the C*-algebra $A(U)$. Another way of writing this
is $\mathcal{F}|_U \subset_\varepsilon \gamma (B)$.
\end{lemma}
\begin{proof}
Let us begin by writing $A(x)$ as  inductive limit of a sequence of
unital Kirchberg algebras with finitely generated K-theory
satisfying the $UCT$, and with unital, injective connecting
$*$-homomorphisms.  This is possible by the classification theorem
of Kirchberg and Phillips (see e.g. \cite[Prop.
8.4.13]{Ror:encyclopedia}). Therefore, there is a unital Kirchberg
subalgebra $B$ with unital inclusion map $\imath: B \subset A(x)$
such that $K_*(B)$ is finitely generated, $B$ satisfies the $UCT$
and $\mathcal{F}|_{\,x} \subset_\varepsilon \imath(B).$  Write
$\mathcal{F} = \{a_1, a_2, ..., a_r\}$ and let $b_1, b_2, ..., b_r
\in B$ such that

\begin{equation}\label{1}\| a_i(x)-\imath(b_i) \| <
\varepsilon,\quad 1 \leq i \leq r.\end{equation}
 Let
$(U_n)_{n \in \mathbb{N}}$ be a decreasing sequence (i.e. $U_{n+1}
\subseteq U_n, n \in \mathbb{N})$, forming a basis system of
compact and open neighborhoods of $x$. In particular $\bigcap_{n
\in \mathbb{N}} U_n = \{x\}$.  Therefore, by Remark \ref{2.6},
$A(x) = \varinjlim (A(U_n), \Phi_n)$, where each $*$-homomorphism
$\Phi_n : A(U_n) \rightarrow A (U_{n+1})$ is the restriction map.
Since $K_*(B)$ is finitely generated and B satisfies the $UCT$
\cite{RosSho:UCT} implies that $KK(B,-)$ is continuous and hence
there exist $n$ and $\alpha\in KK(B, A(U_n))$ such that
$[\pi_n]\alpha = [\imath] \in KK(B, A(x))$, where $\pi_n : A(U_n)
\rightarrow A(x)$ is the restriction map. Now, by Lemma \ref{3.1},
$A(U_n) \cong A(U_n) \otimes \mathcal{O}_\infty$.  Using this, the
fact that $A(U_n)$ has an approximate unit of projections
$(p_n)_{n \in \mathbb{N}}$ (see Remark \ref{2.7} and Remark
\ref{2.2}), and the fact that for every $k\in \mathbb{N}$, every
projection in $M_k \otimes \mathcal{O}_\infty$ is (Murray-von
Neumann) equivalent to a projection in $O_\infty = (e_{11} \otimes
1)(M_k \otimes \mathcal{O}_\infty)(e_{11} \otimes 1)$, it follows
that every projection in $M_k \otimes A(U_n) \cong M_k \otimes
A(U_n) \otimes \mathcal{O}_\infty  \cong \varinjlim M_k \otimes
p_m A(U_n)p_m \otimes \mathcal{O}_\infty$ is equivalent to a
projection in $\varinjlim p_m A(U_n)p_m \otimes \mathcal{O}_\infty
\cong A(U_n)\otimes \mathcal{O}_\infty \cong A(U_n)$. Taking into
account this observation and Kirchberg's theorem \cite[Thm.
8.3.3]{Ror:encyclopedia}, it follows that there is a
$*-$homomorphism $\sigma : B \rightarrow A(U_n)$ with $[\sigma]=
\alpha \in KK (B, A(U_n)).$ From this we obtain that $\pi_n \sigma
: B \rightarrow A(x)$ is a unital $*$-homomorphism which has the
same KK-theory class as the unital $*$-homomorphism $\imath:B
\rightarrow A(x)$.  By the uniqueness part of Kirchberg's theorem,
$\pi_n\sigma$ is approximately unitarily equivalent with $\imath$.
In particular, there is a unitary $u_0 \in A(x)$ such that:
\begin{equation}\label{2} \| u_0 (\pi_n \sigma (b_i)){u_0}^*-\imath(b_i)\| <
\varepsilon, 1\leq i \leq r. \end{equation}
 From \eqref{1} and \eqref{2}  we
obtain:
\begin{equation}\label{3} \| u_0 (\pi_n \sigma (b_i)){u_0}^*-a_i(x) \|< 2
\varepsilon, 1\leq i \leq r. \end{equation} The above inequality
holds in $A(x)$. After increasing $n$ if necessary, we may assume
that $u_0$ lifts to a unitary $u \in A(U_n)$, thus $\pi_n (u)=
u(x)=u_0$.  Using the continuity property of the norm for
continuous fields of C*-algebras, we obtain from (3) after
increasing $n$ if necessary:  \begin{equation}\label{4} \| u(y)
\,\sigma (b_i) (y)\,u(y)^*-a_i(y)\| < 2\varepsilon ,\quad y \in
U_n, \quad 1 \leq i \leq r.
\end{equation}
 If we define $ \gamma : B \rightarrow A(U_n)$ by $\gamma(b)= u \sigma
(b)u^*$ for every $b \in B$ and if we let $ a_i|_{U_n}$ denote the
image of $a_i$ in $A(U_n)$, then \eqref{4} becomes:
\begin{equation}\label{5} \|
\gamma(b_i)-a_i|_{U_n} \| < 2\varepsilon, 1\leq i \leq r
\end{equation}
(since $U_n$ is compact).  Finally, note that \eqref{5} implies
that $\mathcal{F}|_{U_n} \subset_{2\varepsilon} \gamma(B)$.

\end{proof}
\begin{theorem}\label{3.3}
Let $\mathcal{A} = ((A(x))_{x\in X}, \Gamma )$ be a continuous
field of Kirchberg algebras satisfying the $UCT$ over a metrisable
zero dimensional locally compact space $X$.  Let $A$ be the
$C^*$-algebra associated to $\mathcal{A}$. Then $A$ admits an
exhaustive sequence consisting of finite direct sums of unital
Kirchberg algebras satisfying the UCT and having finitely
generated K-theory groups.
\end{theorem}
\begin{proof}
By Remark \ref{2.7} and Remark \ref{2.2} it follows that $A$ has
an approximate unit of projections $ (e_n)_{n \in \mathbb{N}}.$
Since $ A = \overline{\cup_{n \in \mathbb{N}} e_n Ae_n }\cong
\varinjlim e_nAe_n$, it clearly suffices to prove the statement
for each $e_nAe_n$, and thus we may assume -and we shall- that $A$
and all the $A(x)$'s are unital. Fix an arbitrary finite subset
$\mathcal{F} = \{a_1, a_2, ... ,a_m\}$ of $A$ and an arbitrary
$\varepsilon > 0$. For each $1 \leq i \leq m$, define $K_i:=\{x
\in X:\|a_i(x)\|\geq \varepsilon\}$. Then clearly  each $K_i$ is a
compact subset of $X$ and so is $K=\cup_{i=1}^m\,K_i$. Using
Lemma~\ref{3.2} and the compactness of $K$ we find clopen subsets
of $X$ denoted by $U_1,U_2,...,U_k$, unital Kirchberg algebras
$B_1,B_2,...,B_k$ satisfying the UCT with $K_*(B_j)$ finitely
generated, $*$-homomorphisms $\gamma_j:B_j \to A(U_j)$, and
$b_{ij}\in B_j$, $1 \leq i \leq m$, $1 \leq j \leq k$ such that
\begin{equation}\label{6}
K \subseteq U:=\bigcup_{i=1}^m U_i
\end{equation}
and
\begin{equation}\label{7} \| a_i|_{U_j}-\gamma_j(b_{ij}) \| \leq
\varepsilon, \quad 1\leq i \leq m, \quad 1 \leq j \leq
k.\end{equation}
 Let $V_1, V_2,..., V_k$ be mutually disjoint clopen subsets of
 $U_1,
U_2,...,U_k$  such that
\begin{equation*}
\bigcup^k_{j=1} V_j = \bigcup^k_{j=1} U_j.\end{equation*} From
\eqref{7}  we easily get that:
\begin{equation*}\label{10} \| a_i|_{V_j}- r_j\gamma_j(b_{ij})
 \| \leq \varepsilon, 1 \leq i \leq m, 1
\leq j \leq k \end{equation*} where $r_j : A(U_j) \rightarrow
A(V_j)$ is the restriction $*$-homomorphism. Therefore we may
assume that $U_1$, $U_2$,...,$U_k$ are mutually disjoint, and in
particular
$$A(X)\cong A(X\setminus U)\oplus A(U_1)\oplus\cdots \oplus A(U_k).$$
 Let
$C_j=\gamma_j(B_j)\subseteq A(U_j)$ and $C=C_1+\cdots+ C_k\cong
\bigoplus^k_{j=1}C_j$. Since $B_j$ is simple, we have either $C_j
\cong B_j$ or $ C_j=0$ hence $C \in \mathcal{B}_{fg-uct}$.
  Observe
that since the $U_j$'s are clopen, $\chi_{U_j}$ is a continuous
function on $X$ for every $ 1 \leq j \leq k$. Define $$ b_i=
\sum^k_{j=1}\chi_{U_j}\gamma_j(b_{ij}) \in C, \quad 1 \leq i \leq
m.$$  Now, note that \eqref{7} implies that:
\begin{equation}\label{11} \| (a_i - b_i)|_U \| \leq
\varepsilon, 1 \leq i \leq m.\end{equation}
   Now, if $x \in
X\setminus U$ then $b_i(x) = 0$, and  $x \in
 X\setminus K$ and hence, using the definition of
$K$, we get $\| a_i(x)\| < \varepsilon \,( 1 \leq i \leq m).$ This
shows that for $ x \in X\setminus{U}$, we have:
\begin{equation}\label{12} \| a_i(x) - b_i(x) \| = \| a_i(x)\| < \varepsilon, 1
\leq i \leq m\end{equation} In conclusion, \eqref{11}  and
\eqref{12} taken together show that:
$$\| a_i - b_i \| \leq \varepsilon, 1 \leq i \leq m$$
where each $b_i\in C$. Since we already argued  that $C\in
\mathcal{B}_{fg-uct}$, this concludes the proof.
\end{proof}
\begin{corollary}\label{field-sat-uct}
Let $\mathcal{A} = ((A(x))_{x\in X}, \Gamma )$ be a continuous
field of simple separable nuclear C*-algebras  satisfying the
$UCT$ over a metrisable zero dimensional locally compact space
$X$.  Let $A$ be the $C^*$-algebra associated to $\mathcal{A}$.
Then $A$ satisfies the UCT.
\end{corollary}
\begin{proof} Since $\mathcal{O}_\infty$ is
KK-equivalent to $\CCC$, we may replace $A$ by $A \otimes
\mathcal{O}_\infty$.
 By
\cite[Cor. 2.8]{KirWas:bundles} $A \otimes \mathcal{O}_\infty$
satisfies the assumptions of Theorem~\ref{3.3}. By \cite{Dad:uct},
any nuclear separable C*-algebra which admits an exhausting
sequence of C*-subalgebras satisfying the UCT will also satisfy
the UCT.
\end{proof}
 \begin{remark} Let $X$ be a metrisable, zero dimensional, locally compact
space and let $\mathcal{A} =((A(x))_{x \in X},\Gamma)$ be a
continuous field of separable C*-algebras such that each fiber
A(x) admits an exhaustive sequence $(A_n(x))_n$  of simple
semiprojective C*-algebras. Let $A$ be the C*-algebra associated
to $\mathcal{A}$. Then $A$ is the inductive limit of a sequence
$(A_n)_n$, where each $A_n$ is a finite direct sum of C*-algebras
of the form $A_m(x)$ ($m$ and $x$ may vary). The proof of this
statement is similar with the proof of Theorem~\ref{3.3} but it is
much simpler.
\end{remark}
\section{The invariant and basic properties}\label{invar}
We will use K-theory with coefficients. For each $m\geq 2$ let
$W_m$ be the Moore space obtained by attaching a two-cell to the
circle by a degree $m$-map.  Fix a base point $*$ in each of the
spaces $W_m$.
 Let $C(W)$ denote the C*-algebra obtained
by adding a unit to $\bigoplus^\infty_{m=2}  C_0 (W_m\setminus
*).$
 Similarly,  let
 $C(W_M)$ denote the unitalization of
$\bigoplus^M_{m=2}  C_0 (W_m\setminus *)$, where $M$ is an integer
$\geq 2$. Define
$$C=C(\mathbb{T})\otimes C(W)\quad \text{and}\quad C_M=C(\mathbb{T})\otimes
C(W_M).$$ Note that we have natural embeddings $C(\mathbb{T})\otimes
C(W_m)\subset C_M\subset C$,\, ($m \leq M$).

The total K-theory group of a C*-algebra $A$ is given by:
$$ \underline{K} (A)= K_*(A)\oplus \bigoplus^\infty_{m=2} K_* (A;
 \mathbb{Z}/m) \cong K_0(A \otimes C)$$
This group is acted on by the set of  coefficient and Bockstein
operations denoted by $\Lambda$. It is useful to consider the
following direct summand of $\uka$:
 $$\underline{K}(A)_M =
K_*(A)\oplus \bigoplus^M_{m=2}K_*(A; \mathbb{Z}/m) \cong K_0(A
\otimes C_M)$$
\begin{remark}\label{uk-fg} Assume that $A$ is a separable
C*-algebra satisfying the UCT.
 If $M$ annihilates the torsion part
of $K_*(A)$, i.e. $M\, \mathrm{Tors}\, K_*(A)=0$, then the map
$$\mathrm{Hom}_{\Lambda}(\uka,\ukb)\rightarrow \mathrm{Hom}_{\Lambda}(\uka_M,\ukb_M)$$
induced by the restriction map $\uka \to \uka_M$ is bijective for
any $\sigma$-unital C*-algebra $B$, \cite[Cor. 2.11]{DadLor:duke}.
If moreover  $K_*(A)$ is a finitely generated group, then $\uka$
is finitely generated as a $\Lambda$-module, \cite[Prop.
4.13]{DadGong:class-rr0}. More precisely  there are $x_1,...,x_r
\in \uka_M$,  such that for any $x \in \uka$, there exist
$k_1,...,k_r \in \mathbb{Z}$ and $\lambda_1,...,\lambda_r \in
\Lambda$ such that
$$x=k_1\lambda_1(x_1)+\cdots+k_r\lambda_r(x_r).$$
\end{remark}
 If $A$ is a
C*-algebra, we denote by $P(A)$ the Murray-von Neumann abelian
semigroup consisting of equivalence classes of projections in
$A\otimes \mathcal{K}$.

We introduce two homotopy invariants, $\inv (A)$ which will be
used for stable C*-algebras $A$ and $\inv_u(A)$ in the general
case.
 To this purpose we consider
the pair  $(J(A),J(A)^+)$ consisting of a semigroup $J(A)$
together with a distinguished subsemigroup $J(A)^+\subset J(A)$
(called the positive subsemigroup of $J(A)$). This is defined as
follows.
 Consider the map:
$$\rho_0 : P(A\otimes C) \rightarrow P(A \otimes \mathcal{O}_2)$$ induced by
the unital $*-$homomorphism $C \rightarrow \mathcal{O}_2, f
\rightarrow \omega(f)1_{\mathcal{O}_2}$ (where $\omega$ is a fixed
character of the commutative $C^*$-algebra C).  Since the spectrum
of $C$ is path connected, it follows that $\rho_0$ is independent
on the choice of the character $\omega$.  We define: $$J (A):= P(A
\otimes \mathcal{O}_2) \oplus \underline{K} (A)$$ and define
$J(A)^+$ to be the image of the map:
$$\rho: P(A \otimes C) \rightarrow P(A \otimes \mathcal{O}_2) \oplus
\underline{K} (A)= J(A),$$ $$\rho [p]:=\rho_0[p] \oplus
[p]_{\underline{K}(A)}.$$ $J(A)$ is a semigroup with unit $(0,0)$
and is also a $\Lambda$-module (Bockstein and coefficient
operations), where $\Lambda$ acts only on the second component
$\underline{K}(A)$. The pair $(J(A),J(A)^+)$ together with the
action of $\Lambda$ on $J(A)$ is denoted by $\inv (A)$.

In abstract terms, the invariant consists of a
$\mathbb{Z}/2$-graded abelian semigroup $J=J^{(0)}\oplus J^{(1)}$
where $J^{(1)}$ is in fact a graded group acted by $\Lambda$,
together with a distinguished subsemigroup $J^+\subset
J^{(0)}\oplus J^{(1)}$. For a C*-algebra $A$ we have
$J(A)^{(0)}=P(A \otimes \mathcal{O}_2)$  and $J(A)^{(1)}=\uka$.

Note that a $*-$homomorphism $\varphi: A \rightarrow B$ induces a
morphism of the invariant
$$\mathrm{Inv} (\varphi):
\mathrm{Inv}(A) \rightarrow \mathrm{Inv}(B)$$ in the sense that
$\mathrm{Inv} (\varphi)$ is a morphism of graded semigroups which
is $\Lambda$-linear and preserves the positive subsemigroups. We
express these properties by saying that $\mathrm{Inv} (\varphi)$
is positive and $\Lambda$-linear. More precisely:
\begin{definition}\label{4.D}
A map $\alpha : \mathrm{Inv}(A) \rightarrow \mathrm{Inv}(B)$ is
said to be positive and $\Lambda$-linear if it has two components
$$\alpha^{(0)}: P(A \otimes \mathcal{O}_2) \rightarrow P(B \otimes \mathcal{O}_2)$$
with $\alpha^{(0)}$ a unit preserving morphism of semigroups
$(\alpha^{(0)}(0)=0)$ and
$$\alpha^{(1)}: \underline{K}(A) \rightarrow \underline{K} (B)$$
where $\alpha^{(1)}$ is graded morphism of $\Lambda$-modules and
$\alpha (J (A)^+) \subseteq J(B)^+$.
\end{definition}
It is clear that $\inv(A)$ does not distinguish between a
C*-algebra $A$ and its stabilization $A \ot \mathcal{K}$.
Following \cite{Ror:cuntz}, we consider the set $P_u(A)$
consisting of unitary equivalence classes of projections in $A$,
where the unitaries are from the unitalization $\widetilde{A}$ if
$A$ is nonunital. As in \cite{Ror:cuntz}, $P_u(A)$ is equipped
with its family of all finite orthogonal sets.
 A finer
invariant is obtained by enlarging $\inv (A)$ to
$$\inv_u(A)=(\inv(A),P_u(A),\tau),$$ where
$\tau:P_u(A)\to\inv(A)$  is given by the composition of the maps
$P_u(A)\to P(A \ot C)$, $[p]_u\mapsto [p\ot 1_C]$, and $\rho:P(A
\ot C) \to \inv (A).$

A morphism $\alpha:\inv_u(A)\to \inv_u(B)$ consists, in addition
to the components $\alpha^{(i)}$ from Def.~\ref{4.D}, of a map
$\alpha_u:P_u(A)\to P_u(B)$ which maps orthogonal sets to
orthogonal sets and makes the following diagram commutative:
$$\xymatrix{
  {P_u(A)} \ar[d]_{\tau} \ar[r]^{\alpha_u}
                & {P_u(B)} \ar[d]^{\tau}  \\
 \inv (A) \ar[r]_{(\alpha^{(0)},\,\alpha^{(1)})}
                & \inv (B)            }$$

\begin{lemma}\label{4.2}If $A$ is a Kirchberg algebra, then:

$(a)\, P(A \otimes \mathcal{O}_2) = P(\mathcal{O}_2) = \{0, \infty
\}$

$(b)\, J (A)= \{0, \infty \} \oplus \underline{K} (A)$

 $(c)\, J (A)^+ = \{(0,0)\}\cup \{(\infty, x) : x \in
\underline{K}(A)\}$

  $(d)\,$  The map $\rho : P(A \otimes C)
\rightarrow J(A)$ is injective. Moreover $\rho$ is injective even
when $A$ is a finite direct sum of Kirchberg algebras and even
when $A \in \mathcal{L}$.
\end{lemma}
\begin{proof}(a) $A$ is either unital or $A \cong A_0 \otimes
\mathcal{K}$, where $A_0$ is a unital Kirchberg algebra (see
Remark \ref{2.2}).  Since  $P(B) = P(B \otimes \mathcal{K})$ for
every C*-algebra $B$, if follows that we may assume that $A$ is a
unital Kirchberg algebra.  But then, since $A$ is simple,
separable, unital and nuclear, a remarkable result of Kirchberg
(\cite[Thm. 7.1.2]{Ror:encyclopedia}) implies that $A \otimes
\mathcal{O}_2 \cong \mathcal{O}_2$. Therefore $P(A \otimes
\mathcal{O}_2)=P(\mathcal{O}_2)$.  Since any two nonzero
projections in $\mathcal{O}_2 \otimes \mathcal{K}$ are equivalent,
it follows that $P(\mathcal{O}_2)= \{0, \infty \}$. In conclusion,
$P(A \otimes \mathcal{O}_2) = P(\mathcal{O}_2) = \{0, \infty \}.$

(b) This  follows immediately from (a).

(c)  Observe that if $B$ is a simple $C^*$-algebra and $X$ is a
compact, connected space, then every nonzero projection $p$ in $B
\otimes C(X)$ is full. To show that $J(A)^+:=\rho(P(A \otimes C))
\subseteq \{(0,0)\}\cup \{(\infty, x) : x \in \underline{K}(A)\}$,
note that if $\rho_0[p]=0$ for some $p \in P(A\otimes C)$, then
$p=0$. To prove the opposite inclusion, observe first that  $(0,0)
= \rho ([0]) \in J(A)^+$. Now fix an arbitrary $x \in
\underline{K}(A).$ Since $A \otimes O_{\infty} \cong A$, it
follows that $(A \otimes C) \otimes O_{\infty} \cong A \otimes C$
and hence $A \otimes C$ is purely infinite (use \cite[Prop.
4.5]{KR1}).  Let $p$ be an arbitrary nonzero projection of $A$.
Then, by  \cite[Thm. 4.16]{KR1} it follows that $p \otimes 1_C$ is
a properly infinite projection in $A \otimes C$.  Moreover, since
A is simple, it follows that $p \otimes 1_C$ is full in $ A
\otimes C.$  Then, by a result of Cuntz \cite{Cuntz:KOn} (see also
\cite[Lemma 4.15]{BK} it follows that there is a full (and
properly infinite) projection $e$ in $A \otimes C$ such that
$[e]_{\underline{K}(A)}= x \in \underline{K}(A) = K_0 (A \otimes
C).$ Observe also that since $e$ is full in $A \otimes C$, we have
$\rho_0 [e] = \infty$. In conclusion, we have $\rho [e] =
\rho_0[e] \otimes [e]_{\underline{K}(A)} = (\infty, x)$.

(d) First  we consider the case when $A$ is a Kirchberg algebra.
Let $p,q$ be projections in $A \otimes C \otimes \mathcal{K}$ such
that $\rho [p]= \rho [q].$ Since $(A \otimes C \otimes
\mathcal{K})\otimes O_\infty \cong A \otimes C \otimes
\mathcal{K}$ (because $A \otimes O_\infty \cong A)$ and since $A
\otimes C \otimes \mathcal{K}$ has an approximate unit of
projections, it follows (as in the proof of Lemma \ref{3.2}) that
$p$ and $q$ are equivalent with projections in $A\otimes C$.
Hence, we may assume that $ p,q \in A \otimes C$. Since
\mbox{$\rho_0 [p]= \rho_0[q]$} it follows that $p$ and $q$ are
simultaneously zero or nonzero. It suffices to consider the case
when they are both nonzero.  In the proof of (c), we showed that $
A \otimes C$ is purely infinite. Then, by \cite[Theorem 4.16]{KR1}
it follows that $p$ and $q$ (being nonzero) are properly infinite
projections in $A \otimes C$.  Since, as observed above, they are
also full, results of Cuntz \cite{Cuntz:KOn} (see also \cite[Lemma
4.15]{BK}) allow us to conclude from $[p]_{\underline{K}(A)}=
[p]_{K_0(A \otimes C)} =[q]_{K_0(A \otimes C)}=
[q]_{\underline{K}(A)} $ that $p$ and $q$ are equivalent in $A
\otimes C$, i.e. $[p]= [q]$ in $P(A \otimes C).$ Hence $\rho$ is
injective.

In the case $A= \oplus^n_{i=1}A_i$ where each $A_i$ is a Kirchberg
algebra, the fact that the map $\rho: P(A \otimes C) \rightarrow
J(A)$ is injective follows from the fact that $\rho$ is injective
if $A$ is simple (the above case) and the observation that $\rho$
and $P$ are additive with respect to direct sums:
$$\rho = \oplus_{i=1}^n\rho_i : P(A \otimes C)= \oplus^n_{i=1} P(A_i
\otimes C) \rightarrow J(A)= \oplus^n_{i=1} J(A_i).$$


Let $B$ be a C*-algebra with an exhausting sequence $(B_n)$ of
C*-subalgebras. Using functional calculus, one shows that any
partial isometry in $B$ can be approximated by  partial isometries
 in $B_n$'s. Applying this observation to $A\otimes C$ and using the
previous cases we obtain that $\rho$ is also injective for $A\in
\mathcal{L}$.
\end{proof}
\begin{lemma}\label{4.3}
If $A=A_1\oplus\dots\oplus A_m$ with $A_i$
 \kalgs,
then
$$(a) \quad J(A)^+=\{((r_1,x_1),\dots ,(r_m,x_m)): r_i \in P(A_i\otimes \mathcal{O}_2), x_i \in \uk(A_i),
r_i=0\Rightarrow x_i=0\}.$$ In particular this shows that $J(A)^+$
is $\Lambda$-invariant.

$(b)$  If $A$ satisfies the UCT and $K_*(A)$ is \fg, then
$\mathrm{Inv} (A)$ has the following semi-projectivity like
property.
 Let $B=B_1\oplus\dots\oplus B_n$
with $B_i$
 \kalgs. Let $\alpha:J(A)\to J(B)$ be a $\Lambda$-linear morphism
 of graded semigroups. Then there are finitely many elements
 $\{t_1,...,t_r\}$ in $J(A)^+$ with the property that if $\alpha(t_j)\in
 J(B)^+$ for all $1 \leq j \leq r$, then $\alpha(J(A)^+)\subset J(B)^+$.
\end{lemma}
 \begin{proof} The first part follows from Lemma ~\ref{4.2}.  To argue for
 the second part, it suffices to assume that both $A$ and $B$ are
 Kirchberg algebras and $A$ satisfies the UCT.  Since $K_*(A)$ is finitely generated, it
 follows that $ \underline{K}(A)$ is finitely generated as a
 $\Lambda$-module \cite[Prop.
4.13]{DadGong:class-rr0}.
  Let $\{x_1,\dots,x_r\}$ be a list of generators.
  Then $t_j:=(\infty, x_j) \in J(A)^+$
 will satisfy the required property.  Indeed, if $ t \in
 J(A)^+, t \neq (0,0),$
 then $t=(\infty, x)$ for some $x \in \underline{K}(A)$
 and  $x=\sum_j \lambda_jx_j$ for some
 $\lambda_j \in \Lambda$. Thus
  $$ t= (\infty, x) = (\infty,
 \sum_j \lambda_j x_j) = \sum_j \lambda_j (\infty, x_j) = \sum_j
 \lambda_j t_j$$    and $$ \alpha(t) = \sum_j \lambda_j
 \alpha(t_j) \in J(B)^+$$  since $ \alpha(t_j) \in
 J(B)^+$ (by hypothesis) and $ \lambda_j \alpha(t_j) \in
 J(B)^+$, because $J(B)^+$ is
 $\Lambda$-invariant.  Since $\alpha ((0,0))=(0,0) \in
 J(B)^+$, the proof is complete.
 \end{proof}
 \begin{lemma}\label{4.4}
    Let $B $ be a separable C*-algebra with exhaustive sequence $(B_n)$.
    For any $z_1,\dots,z_k\in J(B)^+$, there exist $m$ and
$y_1,\dots,y_k\in J(B_m)^+$  such that map $J(\jmath_m):J(B_m)\to
J(B)$ satisfies $J(\jmath_m)(y_i)= z_i$ .
\end{lemma}
 \begin{proof} This follows by functional calculus
 like in the proof of continuity of $K_0$.
 \end{proof}
 \begin{lemma}\label{4.5}
   If $A$ is a separable C*-algebra with an exhausting sequence $(A_n)$
and all $A_n$ have finitely generated K-theory and satisfy the
UCT, then
   there is a  subsequence $(A_{r(n)})$  of $(A_n)$
   and $\Lambda$-linear morphisms $\nu_n:\uk(A_{r(n)})\rightarrow
   \uk(A_{r(n+1)})$ \st\
   \st\ for each $n$ the diagram
\[
\xymatrix{ \uk(A_{r(n)})\ar[r]^{\nu_n}\ar[rd]_{\uk(\imath_{r(n)})}
&\uk(A_{r(n+1)})\ar[d]^{\uk(\imath_{r(n+1})}\\
&\uk(A)}
\]
is commutative and the induced map $\varinjlim
(\uk(A_{r(n)}),\nu_n)\rightarrow \uk(A)$ is an isomorphism of
$\Lambda$-modules.\end{lemma}
\begin{proof}
    If $A$ is a C*-algebra with $K_*(A)$ \fg\ and satisfies the UCT, then $\uk(A)$ is
    \fg\ as a $\Lambda$-module, actually generated by the group
    generators of $\uk(A)_M$ for some $M$, see Remark~\ref{uk-fg}.
    On the other hand
only finitely elements of $\Lambda$ act on $\uk(A)_M$.
     Hence
    $\uk(A)_M$ is
    a finitely presented $\Lambda$-module. That is,
    one has the group relations among generators and a finite number of relations
    involving finitely many Bockstein operations.
    These observations apply to $A_{r(n)}$.
    Therefore one proves this lemma similarly to the proof of
    Lemma~\ref{4.4}. At each stage one  produces a diagram as in
    the statement which commutes when restricted to the
    $\uk(A_{r(n)})_{M(n)}$, where $M(n)\, \mathrm{Tors} \,K_*(A_{r(n)})=0$.
    Finally one applies Remark~\ref{uk-fg}.
\end{proof}

 \begin{lemma}\label{4.6}
Let $A$ be a C*-algebra which admits an exhaustive sequence
$(A_n)$ such that each $A_n$ is a finite direct sum of \kalgs\
satisfying the UCT and  with finitely generated K-theory. Then
after passing to a subsequence of $(A_n)$ (if necessary) there is
an inductive system
\begin{equation*}
    \cdots (J(A_n),J(A_n)^+)\stackrel{\nu_n}\rightarrow(J(A_{n+1}),J(A_{n+1})^+)\cdots
\end{equation*}
\st\ for each $n$, $\nu_n$ is positive and $\Lambda$-linear, the
diagram
\[
\xymatrix{
J(A_n)\ar[r]^{\nu_n}\ar[rd]_{J(\imath_n)}&J(A_{n+1})\ar[d]^{J(\imath_{n+1})}\\
&J(A)}
\]
is commutative and the induced map $\varinjlim
((J(A_n),J(A_n)^+),\nu_n)\rightarrow(J(A),J(A)^+)$ is a positive
$\Lambda$-linear isomorphism. One can replace $(J(-),J(-)^+)$ by the
alternate notation $\mathrm{Inv}(-)$ everywhere in the statement of
this lemma, since the action of $\Lambda$ was considered.
\end{lemma}
\begin{proof}
This follows by putting together Lemma~\ref{4.3}, Lemma~\ref{4.4},
Lemma~\ref{4.5}. The property given in the  part $(b)$ of
Lemma~\ref{4.3} is crucial since it allows to insure positivity of
$\Lambda$-linear maps by obtaining it only for finitely many
elements.
\end{proof}
\begin{lemma}\label{4.7}  Let $A$ be a $C^*$-algebra which admits an exhaustive
sequence ($A_n)$ such that each $A_n$ is a finite direct sum of
Kirchberg algebras satisfying the UCT and  with finitely generated
K-theory. Let $B$ be a $C^*$-algebra which admits an exhaustive
sequence $(B_n)$ such that each $B_n$ is a finite direct sum of
Kirchberg algebras. If $\alpha: \mathrm{Inv}(A) \rightarrow
\mathrm{Inv}(B)$ is a positive $\Lambda$-linear map, for any $n$,
there exist $m=m(n)$ and a positive $\Lambda$-linear map
$\alpha_n: \mathrm{Inv}(A_n) \rightarrow \mathrm{Inv} (B_m)$ such
that the diagram
$$\xymatrix{
  {\inv (A_n)} \ar[d]_{\alpha_n} \ar[r]^{\inv (\imath_n)}
                & \inv (A) \ar[d]^{\alpha}  \\
 \inv (B_m) \ar[r]_{\inv(\jmath_m)}
                & \inv (B)            }$$
   is commutative.
\end{lemma}
\begin{proof}  This is just a repetition of the proof of
Lemma \ref{4.6}. The crucial property of $\mathrm{Inv}(-)$ is that
it is  "finitely presented" when applied to finite direct sums of
Kirchberg algebras with finitely generated K-theory.
\end{proof}
\begin{remark} If $A$ is a unital Kirchberg algebra, then
$$P_u(A)\cong \{[0]_u, [1]_u\}\sqcup K_0(A),$$
by \cite{Cuntz:KOn}, \cite{Ror:cuntz}. It is then clear that by
using similar arguments one shows that $\inv_u(-)$ satisfies
statements analog to Lemmas~\ref{4.6} and \ref{4.7}.
\end{remark}

\section{Classification results}\label{section-class}

    We begin by recalling some terminology and definitions from
    \cite{DadEil:class}
    which will be important in what follows.

    \begin{definition}\label{5.1} (\cite{DadEil:class}).  Let $A$ be a
    C*-algebra.  A $\underline{K}$-triple $(\mathcal{P},
    \mathcal{G}, \delta)$ consists of a finite subset $\pset$ of projections,
    $$\mathcal{P}\subseteq \bigcup_{m\geq 1} Proj
    (A\otimes C(\mathbb{T}) \otimes C(W_m) \otimes
    \mathcal{K})$$
    where the $W_m$'s are the Moore spaces of order $m$,
    and a finite subset $\mathcal{G} \subseteq A$, and  $\delta > 0$ chosen such
    that whenever $\varphi$ is a completely positive contraction
    which is $\delta$-multiplicative on $\mathcal{G}$, then
    $e=(\varphi \otimes id) (p)$ is almost a projection
    in the sense that $\|e^2-e\|<1/4$ and in particular
    $$ 1/2 \notin sp((\varphi \otimes id) (p))$$
    for each $p \in \mathcal{P}$, where $id$ is the identity of
    $C(\mathbb{T}) \otimes C(W_m) \otimes \mathcal{K}$ for suitable $m$.
    \end{definition}
\begin{definition}\label{5.2}(\cite[Def. 3.9]{DadEil:class}  Let $(\mathcal{P},
\mathcal{G}, \delta)$ be a $\underline{K}$-triple, and assume that
$\varphi:A \rightarrow B$ is a completely positive contraction
which is $\delta$-multiplicative on $\mathcal{G}$.  We define
$\varphi_{\sharp} (p)=[\chi_0 (\varphi \otimes id) (p)]$ where
$\chi_0 : [0,1] \setminus \{1/2\} \rightarrow [0,1]$ is $0$ on
$[0, {1/2})$ and $1$ on $({1/2},1]$. It is clear that if
 $(\mathcal{P}, \mathcal{G}, \delta)$ is a $\uka$-triple one also has a
 natural map
$\varphi_{\sharp}: \mathcal{P} \rightarrow P(B \otimes C)$
whenever $\varphi : A \rightarrow B$ is a completely positive
contraction which is $\delta$-multiplicative on $\mathcal{G}$.
Hence, one can define in this way $\varphi_{\sharp} : \mathcal{P}
\rightarrow \inv(B)$ whenever $(\mathcal{P}, \mathcal{G}, \delta)$
is a $\underline{K}(A)$ -triple and $\varphi : A \rightarrow B$ is
a completely positive contraction which is $\delta$-multiplicative
on $\mathcal{G}$. Similarly one defines a map $\varphi_{\sharp}
:\pset \to P_u(B)$.

\end{definition}
We need the following result.
\begin{theorem}\label{5.3} Assume that A is a unital Kirchberg algebra which satisfies
the UCT.  For any finite set $\mathcal{F} \subset A$ and
$\varepsilon > 0$, there is a $\underline{K}(A)$-triple
$(\mathcal{P}, \mathcal{G}, \delta)$ with the following property.
For any unital Kirchberg algebra $B$ and any $(\mathcal{G},
\delta)$-multiplicative completely positive contractions
$\varphi,\psi:A \rightarrow B$ with
$\varphi_{\sharp}(p)=\psi_{\sharp}(p) \in \underline{K}(B)$ for
all $p \in \mathcal{P}$, there is an unitary $u \in B$  such that
$\| u \varphi(a)u^*-\psi(a)\|<\varepsilon$ for all $a\in
\mathcal{F}$.
\end{theorem}
\begin{proof}
This is \cite[Thm. 6.20]{DadEil:class}. In addition we may work
with a set of projections $\pset\subset A \otimes C$ rather than
$\pset\subset A \otimes C\otimes \mathcal{K}$.
\end{proof}
\begin{corollary}\label{5.4}  Assume that $A$ is a finite direct sum of
Kirchberg algebras that satisfies the UCT.  For any finite set
$\mathcal{F}\subset A$ and $\varepsilon > 0$, there is a
$\uk(A)$-triple $(\mathcal{P}, \mathcal{G}, \delta)$ with the
following property. For any  finite direct sum of stable Kirchberg
algebras $B$ and any $(\mathcal{G}, \delta)$-multiplicative
completely positive contractions $\varphi, \psi: A \rightarrow B$
with $\varphi_{\sharp}(p) = \psi_{\sharp}(p) \in \inv (B)$ for all
$p\in \mathcal{P}$, there is a unitary $u\in \widetilde{B}$ (the
unitalization of $B$) such that $\| u\varphi(a)u^{*}-\psi(a)\| <
\varepsilon$ for all $a\in \mathcal{F}$.
\end{corollary}
\begin{proof}First observe that we may assume that $B$ is simple
(otherwise we compose with the projection onto each direct
summand). Now, since every Kirchberg algebra has an approximate
unit of projections, we may assume that $A= {\oplus}^{k}_{i=1}A_i$
where each $A_i$ has a unit denoted by $e_i$.  We choose
$\mathcal{P}$ such that $e_i \in \mathcal{P}, 1\leq i\leq k$.
After a small perturbation we may assume that
$(\varphi(e_i))^{k}_{i=1}$ and $(\psi(e_i))^{k}_{i=1}$ are finite
sequences consisting each of mutually orthogonal projections.
Since $\varphi(e_i)$ and $\psi(e_{i})$  have the same class in
$P(B \otimes \mathcal{O}_{2})$ they are simultaneously zero or
nonzero.  Since they have the same class in $\underline{K}(B)$, in
particular they have the same class in $K_{0}(B)$.  But $B$ is
purely infinite and simple. By results of Cuntz \cite{Cuntz:KOn},
it follows that $\varphi (e_{i})$ is equivalent to $\psi (e_{i})$,
$1\leq i\leq k$. Hence, after conjugating $\varphi$ with a
suitable partial isometry $v\in B$, we may assume that $\varphi
(e_{i})=\psi (e_{i}), 1\leq i\leq k$. Since $B$ is stable, this
partial isometry extends to a unitary $u \in \widetilde{B}$.
 Next  we apply
Theorem \ref{5.3} for the restriction of the two given maps to
$A_i \rightarrow \varphi (e_{i}) B\varphi (e_{i}),1\leq i\leq k$.
Thus
 we obtain a unitary $w\in eBe$, where $e=\varphi(e_1)+\cdots+\varphi(e_k)$, with
$\| w\varphi(a)w^{*}-\psi(a)\| < \varepsilon$ for all $a\in
\mathcal{F}$. Finally we choose $u=1-e+w\in U(\widetilde{B})$.
\end{proof}


 The following is our "uniqueness" result
\begin{proposition}\label{5.6} Assume that $A\in \mathcal{L}_{uct}$.  For any
finite set $\mathcal{F}\subset A$ and $\varepsilon > 0$, there is
a $\uk(A)$-triple $(\mathcal{P},\mathcal{G},\delta)$ with the
following property.  For any stable $B\in \mathcal{L}$ and any
$(\mathcal{G}, \delta)$- multiplicative completely positive
contractions $\varphi, \psi : A\rightarrow B$ with
$\varphi_{\sharp}(p) = \psi_{\sharp}(p)\in \inv(B)$ for all $p\in
\mathcal{P}$, there is a unitary $u\in \widetilde{B}$ such that
$\| u\varphi(a)u^{*}-\psi(a)\| < \varepsilon$ for all $a\in
\mathcal{F}$. If $B$ is not assumed to be stable, the conclusion
remains valid if we assume that $\varphi_{\sharp}(p)\in P_u(B)$
for all $p\in \mathcal{P}\cap A$.
\end{proposition}
\begin{proof}Let $(A_{n})$ be an exhaustive
sequence for $A$ with $A_{n} \in \mathcal{B}_{uct}$. It suffices
to prove the statement for the restrictions of $\varphi$ and
$\psi$ to $A_n$. Thus we reduced the proof to the case when $A \in
\mathcal{B}_{uct}$.

 Using Lemma \ref{5.5} we
find an exhausting sequence $(B_{n})$ for $B$ with $B_{n}\in
\mathcal{B}$ and a sequence of completely positive contractions
$\pi_{n}:B\rightarrow B_{n}$ which is asymptotically
multiplicative and such that:
\begin{equation}\label{18}
lim_{n\rightarrow \infty}\| \jmath_{n}\pi_{n}(b)-b\| = 0,
\forall\,b\in B.\end{equation}
  Define for every $n\in \mathbb{N},
\varphi_{n},\psi_n:A\rightarrow B_{n}$ by
$\varphi_{n}:=\pi_{n}\varphi$, $\psi_{n}:=\pi_{n}\psi$. Then
\eqref{18} implies that for $\forall\,a\in A$, we have:
$$lim_{n\rightarrow \infty} \|
\jmath_{n}\varphi_{n}(a)-\varphi (a)\| = 0,$$ A similar property
is satisfied by $\psi_n$.
  Replacing $\varphi$ with $\varphi_{n}$ and
similarly, $\psi$ with $\psi_{n}$, for large $n$, we reduce the
proof of the proposition to Corollary ~\ref{5.4}.
\end{proof}
\begin{lemma}\label{5.7}  Let
$A\in \mathcal{B}_{fg-uct}$ and let $B$ be an finite direct sum of
Kirchberg algebras.  Then any positive and $\Lambda$-linear
morphism $\alpha :\inv(A)\rightarrow \inv(B)$  lifts to a
$*$-homomorphism $\varphi :A\rightarrow B$.

\end{lemma}
\begin{proof}We may suppose that $B$ is a Kirchberg algebra.  Write
$A=\bigoplus_{i=1}^{k}A_{i}$ with $A_{i}\in \AAA_{fg-uct}$. Then,
clearly:
$$\inv(A)= \oplus_{i=1}^{k}\inv(A_{i})$$
and correspondingly, $\alpha_{i}$ is the restriction of $\alpha$
to $\inv(A_{i})$.  All we have to do is to lift each $\alpha_{i}$
to a *-homomorphism $\varphi_{i}:A_{i}\rightarrow B$.  Then, since
$B$ is purely infinite and simple, using \cite{Cuntz:KOn} we may
assume, after conjugating each $\varphi_{i}$ with a suitable
partial isometry in $B$, that
$\varphi_{i}(1_{A_i})\varphi_{j}(1_{A_j})=0$ whenever $i\neq j$.
Finally, we shall define $\varphi :A\rightarrow B$ as the direct
sum of the $\varphi_{i}$'s.
    Now, as explained earlier:
    $$J(A_{i})^{+}=\{(0,0)\}\cup \{\infty \oplus
    \underline{K}(A_i)\}\subset P(A_i \otimes \mathcal{O}_2) \oplus
    \underline{K}(A_i)$$
    and:
    $$J(B)^+ =\{(0,0)\}\cup \{\infty \oplus \underline{K}
    (B)\}\subset P(B\otimes \mathcal{O}_2) \oplus
    \underline{K}(B).$$
    The map $\alpha_i$ has components $\alpha_i^{(0)}: P(A_i \otimes
    \mathcal{O}_2)\rightarrow P(B\otimes \mathcal{O}_2 )$ and
    $\alpha_i^{(1)}:\underline{K}(A_i)\rightarrow \underline{K}(B),$
    $$\alpha_i=\left(%
\begin{array}{cc}
  \alpha_i^{(0)} & 0\\
  0 & \alpha_i^{(1)}\\
\end{array}%
\right).
$$
    Therefore the condition $\alpha_i (J(A_i)^+)\subseteq J(B)^+$
    shows that $\alpha_i^{(1)}=0$ whenever
    $\alpha_i^{(0)}(\infty)=0$.  If $\alpha_i^{(0)}(\infty)=0$, then
    $\alpha_i^{(0)}=0$ and we set $\varphi_i =0, 1\leq i\leq n$.  If
    $\alpha_i^{(0)}=\infty$, we apply the UMCT of \cite{DadLor:duke} to lift
    $\alpha_i^{(1)}$ to an element $\beta_i$ of $KK(A_i, B)$.  Since
    $K_*(A_i)$ is finitely generated, $A_i$ satisfies the $UCT$
    and $B$ has a countable approximate unit of projections, by using the continuity
    of $KK(A_i,-)$
    we may assume that $B$ is unital.  Now, using \cite[Thm. 8.3.3]{Ror:encyclopedia},
     we lift $\beta_i$ to a full $*$-homomorphism
    $\varphi' _i:A_i \rightarrow B\otimes \mathcal{K}$.  But, notice
    that every projection in a matrix algebra over $B$ is equivalent
    to a projection in  $B\subset B\otimes \mathcal{K}$ (since
    $B\otimes \OOO_{\infty}\cong B$ and $B$ has an approximate unit of
    projections).  Hence, there is a partial isometry $v_i \in
    B\otimes \mathcal{K}$ such that $v_i^*v_i=\varphi_i'(1_{A_i})$
    and $v_iv_i^* \leq 1_B$.  Define $\varphi_i : A_i\rightarrow B$
    by $\varphi_i : = v_i\varphi_i' v_i^{*}$.  Then, clearly,
    $[\beta_i]$ lifts to $\varphi_i$ and hence $\alpha_i^{(1)}$
    lifts to $\varphi_i$.  It is clear that $\varphi_i:A_i \rightarrow
    B$ is injective ($\varphi_i$ is nonzero and $A_i$ is unital) and
    hence $\varphi_i$ also lifts $\alpha_i^{(0)}$.

\end{proof}
The following is our "existence" result:
\begin{proposition} \label{5.8} Assume that $A\in \mathcal{L}_{uct}$.  For
any $\uk(A)$-triple $(\mathcal{P}, \mathcal{G}, \delta)$, any
$B\in \mathcal{L}$ and any positive $\Lambda$-linear maps
$\alpha:\inv(A)\rightarrow \inv(B)$, there is a
$(\gset,\delta)$-multiplicative completely positive contraction
$\varphi:A \to B$ such that $\varphi_{\sharp}(p) = \alpha ([p])$
for all $p\in \mathcal{P}$.
\end{proposition}
\begin{proof}
Note that $\AAA_{uct} \subseteq \mathcal{L}_{fg-uct}$ (see e.g.
\cite[Prop. 8.4.13]{Ror:encyclopedia}) hence
$\mathcal{L}_{uct}=\mathcal{L}_{fg-uct}$. Let $(A_n)$ be an
exhaustive sequence for $A$, such that $A_n \in
\mathcal{B}_{fg-uct}$ for each $n$.  Passing to a
  subsequence, we may assume that $(A_n)$  satisfies \eqref{exhaustive}.
  By Lemma ~\ref{5.5}, there is a sequence of
completely positive contractions $\mu_n: A\rightarrow A_n$ which
is asymptotically multiplicative and such that :
\begin{equation}\label{19}lim_{n\rightarrow \infty}
\parallel \imath_n\mu_n(a)-a\parallel = 0, \forall\,a\in A.\end{equation}
Applying Lemma \ref{4.7}, we obtain a commutative diagram:
$$\xymatrix{
  {\inv (A_n)} \ar[d]_{\alpha_n} \ar[r]^{\inv (\imath_n)}
                & \inv (A) \ar[d]^{\alpha}  \\
 \inv (B_m) \ar[r]_{\inv(\jmath_m)}
                & \inv (B)            }$$
 where
$\alpha_n :\inv(A_n)\rightarrow \inv(B_m)$ is positive and
$\Lambda$-linear, for some $m=m(n)$.  Next we lift $\alpha_n$ to a
$*$-homomorphism $\varphi'_{n,m}: A_n \rightarrow B_m$ using Lemma
\ref{5.7}.  Finally we set $\varphi_{n,m}: =
\jmath_m\varphi'_{n,m}\mu_n : A \rightarrow B$. Using \eqref{19},
one verifies that if $n$ is large enough, $\varphi:=\varphi_{n,m}$
will satisfy the conclusion of the theorem.
\end{proof}
\begin{theorem} \label{5.9} Let $A, B$ be stable $C^*$-algebras which admit
exhaustive sequences consisting of finite direct sums of Kirchberg
algebras satisfying the UCT. Then $A$ is isomorphic to $B$ if and
only if $\inv(A) \cong \inv(B)$ as $\Lambda$-modules. Moreover
$$Hom(A,B)/\overline{Inn}(B) \cong
Hom_{\Lambda}(\inv(A) ,\inv(B))$$ if $A \in \mathcal{L}_{uct}$ and
$B \in \mathcal{L}$.
\end{theorem}
\begin{proof}  We will only prove the first part, as the second part is similar.
Fix a positive $\Lambda$-linear isomorphism $\alpha : \inv(A)
\rightarrow \inv(B)$.  We may apply the "existence" result
Proposition \ref{5.8} to get completely positive contractions
$\varphi_i : A\rightarrow B$ and $\psi_i :B\rightarrow A$ which
are increasingly multiplicative on larger and larger sets, and
induce $\alpha$ and $\alpha^{-1}$, respectively, on larger and
larger subsets of $\inv(A)$ and $\inv(B)$. Arranging this
appropriately, we may conclude by our "uniqueness" result
Proposition \ref{5.6} that there are unitaries $u_n$ and $v_n$
making
$$ \xymatrix{
A \ar[rr]^{\Ad{u_1}} \ar[dr]_{\varphi_{1}}&
& A \ar[rr]^{\Ad{u_2}} \ar[dr]_{\varphi_{2}}& & A \ar[dr] \ar[r] & ... \\
& B \ar[rr]_{\Ad{v_1}} \ar[ur]^{\psi_{1}}& & B\ar[ur]^{\psi_{2}}
\ar[rr]_{\Ad{v_2}}& & B \ar[r]& ...}
$$
an approximate intertwining in the sense of Elliott
\cite{Ell:rrzeroI}.  Hence
$$A \cong \varinjlim (A, \Ad{u_n}) \cong \varinjlim (B, \Ad{v_n})
\cong B.$$
\end{proof}
\begin{theorem} \label{5.99}
Let $A$, $B$ be stable  C*-algebras associated to  continuous
fields of Kirchberg algebras satisfying the $UCT$ over zero
dimensional metrisable locally compact spaces. Then $A$ is
isomorphic to $B$ if and only if $\inv(A)\cong \inv(B)$.
\end{theorem}
\begin{proof}
    This follows from Theorems \ref{3.3} and \ref{5.9}.
\end{proof}
\begin{corollary}\label{compare}
Let $A$, $B$ be a
 separable nuclear  C*-algebras with
  zero dimensional  Hausdorff primitive spectra.
Assume that all simple quotients of $A$ and $B$ satisfy the UCT.
The  following assertions are equivalent:

$(a)$ \mbox{$A\otimes\mathcal{O}_\infty \otimes\mathcal{K}\cong
B\otimes\mathcal{O}_\infty \otimes\mathcal{K}$};

$(b)$ $\inv(A)\cong \inv(B)$;

$(c)$ $A$ is $KK_X$-equivalent to $B$.
\end{corollary}
\begin{proof} The equivalence (a) $\Leftrightarrow$ (c) is due to Kirchberg and
holds in much more generality \cite{Kir:Michael}. The novelty here
is (a) $\Leftrightarrow$ (b) for which we give a direct proof and
hence obtain an algebraic criterion for when (c) happens to hold.
Let $X:=\mathrm{Prim}(A)$.
 By
Fell's theorem \cite{Fell}, $A$ is isomorphic to a the C*-algebra
associated to a continuous field of C*-algebras
$\mathcal{A}=((A(x))_{x\in X},\Gamma)$ with each $A(x)$ simple.
Then $A \otimes \mathcal{O}_\infty\otimes \mathcal{K}$ is a
C*-algebra satisfying the assumptions of Theorem~\ref{5.99} and
$B\otimes \mathcal{O}_\infty\otimes \mathcal{K}$ has a similar
property, see \cite[Cor. 2.8]{KirWas:bundles}. After this
preliminary discussion note that
 (a) $\Leftrightarrow$ (b) follows from Theorem~\ref{5.99}.
\end{proof}

\begin{corollary} \label{5.10} Let $A$ be a stable $C^*$-algebra which
admits exhaustive sequence consisting of finite direct sums of
Kirchberg algebras satisfying the UCT. Then $A \cong \varinjlim
A_n$, where each $A_n$ is a finite direct sum of Kirchberg
algebras satisfying the UCT and having finitely generated
K-theory. In particular this applies to the stabilization of the
$C^*$-algebra associated to a continuous field of Kirchberg
algebras satisfying the UCT, over a metrisable zero dimensional
locally compact space.
\end{corollary}
\begin{proof} Since $A \in \mathcal{L}_{uct} = \mathcal{L}_{fg-uct}$, it
follows that there is an exhausting sequence $(A_n)$ for $A$, with
$A_n \in \mathcal{B}_{fg-uct}$. By Lemma \ref{4.6}, after passing
to a  subsequence of $(A_n)$ if necessary, we can arrange that
$\inv(A) \cong \varinjlim (\inv(A_n), \nu_n)$, for some positive
$\Lambda$-linear maps $\nu_n:\inv(A_n) \rightarrow\inv( A_{n+1})$
$ (n \in \mathbb{N})$.  Now, using Lemma \ref{5.7}, we lift each
$\nu_n$ to a $*$-homomorphism $\Phi_n : A_n \rightarrow A_{n+1}$.
Define $B': = \varinjlim  (A_n, \Phi_n)$ and  $B: = B'\otimes
\mathcal{K}$. Since $B = \varinjlim A_n\otimes M_n$ and
$A_n\otimes M_n\in \mathcal{B}_{fg}$ (since $A_n \in
\mathcal{B}_{fg}$) for each $n$, it follows obviously that $B \in
\mathcal{L}_{uct}$. Hence $A, B$ are stable $C^*$-algebras in
$\mathcal{L}_{uct}$ and
$$\mathrm{Inv}(B) = \mathrm{Inv}(B'\otimes \mathcal{K}) \cong
\mathrm{Inv}(B') \cong \varinjlim (\,\mathrm{Inv}(A_n),
\mathrm{Inv}(\Phi_n)) = \varinjlim (\,\mathrm{Inv}(A_n), \nu_n)
\cong \mathrm{Inv}(A).$$ Then Theorem \ref{5.9} implies that $A
\cong B \cong \varinjlim (A_n, \Phi_n)$, which ends the proof,
since each $A_n \in \mathcal{B}_{fg-uct}$.
\end{proof}
\begin{corollary}\label{5.11}  Let $A, B$ be stable $C^*$-algebras which admit
exhaustive sequences consisting of finite direct sums of Kirchberg
algebras satisfying the UCT.  Then, the following assertions are
equivalent:

$(a)\, A$ and $B$ are $*$-isomorphic;

$(b) \,A$ and $B$ are shape equivalent;

$(c)\, A$ and $B$ are homotopy equivalent.
\end{corollary}
\begin{proof}Since the implications $(a) \Rightarrow (b)$, $(a)
\Rightarrow (c)$ and $(c) \Rightarrow (b)$ are trivially true, to
prove the corollary it is enough to prove that $(b) \Rightarrow
(a)$ is true. $\inv(-)$ is a homotopy invariant continuous
functor. Therefore if $A$ and $B$ are shape equivalent, then we
clearly have $\mathrm{Inv}(A) \cong \mathrm{Inv}(B)$, which by
Theorem \ref{5.9} implies that $A$ and $B$ are $*$-isomorphic.
\end{proof}
\begin{remark} Theorems \ref{5.9} and \ref{5.99} remain true for nonstable
C*-algebras, provided that one replaces $\inv(-)$ by $\inv_u(-)$ in
their statements. The proof is essentially the same, except for
small changes as in \cite{Ror:cuntz}. Consequently, the stability
assumptions can be dropped from Corollaries ~\ref{5.10} and
\ref{5.11}. For Corollary ~\ref{5.10}, which is probably the more
interesting statement of the two, one can also verify our claim, at
least in the unital case, as follows. Assume that $A\in \LLL_{uct}$
is unital. Then, by the stable case, we can write $A \otimes
\mathcal{K}$ as the closure of an increasing sequence $(A_n)$ of
C*-algebras in $\mathcal{B}_{fg-uct}$. Letting $p=1_A\otimes
e_{11}$, we may assume that $\|p-q\|<1$ for some projection $q \in
A_1$. Let $v \in A \otimes \mathcal{K}$ be a partial isometry with
$v^*v=p$ and $vv^*=q$. Thus $\mathrm{Ad(v)}: p(A\otimes
\mathcal{K})p\to q(A \otimes \mathcal{K})q$  is an isomorphism. Then
$$A \cong p(A\otimes \mathcal{K})p\cong q(A \otimes
\mathcal{K})q=\overline{\cup_n qA_nq}$$ and
$qA_nq\in\mathcal{B}_{fg-uct}$.
\end{remark}
\emph{Acknowledgements}. The first author is partially supported
by NSF grant DMS-0200601 and the second author is partially
supported by NSF grant DMS-0101060.  This material is based upon
work supported by, or in part by, the U.S. Army Research Office
under grant number DAAD19-00-1-0152 (for the second author).

\bibliographystyle{amsplain}

\providecommand{\bysame}{\leavevmode\hbox
to3em{\hrulefill}\thinspace}
\providecommand{\MR}{\relax\ifhmode\unskip\space\fi MR }
\providecommand{\MRhref}[2]{%
  \href{http://www.ams.org/mathscinet-getitem?mr=#1}{#2}
} \providecommand{\href}[2]{#2}

 \end{document}